\documentclass[12pt]{amsart} 

\usepackage{amsmath} 
\usepackage{amssymb} 
\usepackage{amsthm}
\usepackage{amscd} 
\usepackage[all]{xy} 

\def\F{{\mathbb F}}

\def\Q{{\mathbb Q}}
\def\R{{\mathbb R}}
\def\Z{{\mathbb Z}}

\def\O{{\mathcal O}}
\def\I{{\mathcal I}}
\def\J{{\mathcal J}}
\def\a{{\mathfrak a}}
\def\b{{\mathfrak b}}
\def\m{{\mathfrak m}}

\def\Ann{\text{\rm Ann}}

\def\Hom{\text{\rm Hom}}

\def\Spec{\text{\rm Spec}\,}
\def\Supp{\text{\rm Supp}\,}
\newcommand{\Div}{\mathrm{div}}
\renewcommand{\labelenumi}{\textup{(\arabic{enumi})}}
\numberwithin{equation}{section}

\theoremstyle{plain}
\newtheorem{thm}{Theorem}[section]
\newtheorem{mainthm}{Theorem}
\newtheorem{cor}[thm]{Corollary}
\newtheorem{prop}[thm]{Proposition}
\newtheorem{defthm}[thm]{Definition-Theorem} 
\newtheorem{lem}[thm]{Lemma}
\newtheorem{conj}[thm]{Conjecture}
\theoremstyle{definition} 
\newtheorem{defn}[thm]{Definition}
\newtheorem{expl}[thm]{Example} 
\theoremstyle{remark}
\newtheorem{rem}[thm]{Remark}
\newtheorem*{claim}{Claim}
\newtheorem*{acknowledgement}{Acknowledgement}

\title{F-singularities of pairs and Inversion of Adjunction of arbitrary codimension}
\author{Shunsuke Takagi}
\address{Graduate School of Mathematical Sciences, University of Tokyo, 3-8-1, Komaba, Meguro, Tokyo 153-8914, Japan}
\email{stakagi@ms.u-tokyo.ac.jp}
\begin{document}
\begin{abstract}
We generalize the notions of F-regular and F-pure rings to pairs $(R,\a^t)$ of rings $R$ and ideals $\a \subset R$ with real exponent $t > 0$, and investigate these properties. These ``F-singularities of pairs'' correspond to singularities of pairs of arbitrary codimension in birational geometry.
Via this correspondence, we prove a sort of Inversion of Adjunction of arbitrary codimension, which states that for a pair $(X,Y)$ of a smooth variety $X$ and a closed subscheme $Y \subsetneq X$, if the restriction $(Z, Y|_Z)$ to a normal $\Q$-Gorenstein closed subvariety $Z \subsetneq X$ is klt (resp. lc), then the pair $(X,Y+Z)$ is plt (resp. lc) near $Z$.
\end{abstract}

\maketitle
\markboth{SHUNSUKE TAKAGI}{INVERSION OF ADJUNCTION OF ARBITRARY CODIMENSION}

\section*{Introduction}
For the classification theory of higher dimensional algebraic varieties, it is important to study the singularities of a pair $(X,Y)$ of a normal variety $X$ and a closed subscheme $Y$ of $X$.
The purpose of this paper is to investigate local properties of the pair $(X,Y)$, especially Inversion of Adjunction of arbitrary codimension, using the theory of tight closure.

Consider a pair $(X, Y)$ of a non-singular variety $X$ over a field of characteristic zero and a formal combination $Y=\sum_{i=1}^kt_iY_i$ where $Y_i \subsetneq X$ are closed subschemes and $t_i >0$ are real numbers. Let $Z \subsetneq X$ be a normal $\Q$-Gorenstein closed subvariety of codimension $r>0$ such that $Z \not\subset \cup_{i=1}^k Y_i$. 
In case $Z$ is a locally complete intersection variety, a result of Koll\'{a}r \cite{K+} and Shokurov \cite{Sh} says that the pair $(Z,Y|_Z)$ is klt if and only if the pair $(X,Y+rZ)$ is plt near $Z$. 
The lc case, which states $(Z,Y|_Z)$ is lc if and only if $(X,Y+rZ)$ is lc near $Z$, was proved 
by Ambro \cite{Am2} when $X$ is an affine space and $Z$ is a non-degenerate hypersurface, 
by Ein, Musta{\c{t}}{\v{a}} and Yasuda \cite{EMY} via the theory of jet schemes and motivic integration when $Z$ is a divisor, and recently by Ein and Musta{\c{t}}{\v{a}} \cite{EM} when $Z$ is a locally complete intersection variety.
However, in the case where $Z$ is not locally complete intersection, $(X,Y+rZ)$ is not necessarily plt (resp. lc) near $Z$, even if $(Z,Y|_Z)$ is klt (resp. lc) (cf. Example \ref{example}).
In this paper, when $Z$ is not necessarily locally complete intersection, assuming that $(Z,Y|_Z)$ is klt (resp. lc), we give a lower bound for the lc-threshold of $Z$ with respect to $(X,Y)$ (we refer to such statement as ``a sort of Inversion of Adjunction of arbitrary codimension for klt (resp. lc) pairs'' in this paper).

The notion of tight closure is a powerful tool in commutative algebra introduced by Hochster and Huneke \cite{HH1} about fifteen years ago.
There is a strong connection between the singularities arising in birational geometry and the singularities obtained from the theory of tight closure. 
For example, for a $\Q$-Gorenstein normal local ring $R$ of characteristic zero, $\Spec R$ has only log terminal singularities if and only if $R$ is of F-regular type, and if $R$ is of dense F-pure type, then $\Spec R$ has only log canonical singularities (see \cite{Ha2}, \cite{HW}, \cite{MS} and \cite{Sm1}). 
As an analog of singularities of pairs in birational geometry, we generalize the notions of F-regular and F-pure rings to those for pairs $(R, \a^t)$ of rings $R$ and ideals $\a \subset R$ with real exponent $t>0$, and even more, those for several ideals (Definition \ref{ideal}). We look into those properties, and we prove a generalization of ``F-inversion of Adjunction'' \cite[Theorem 4.9]{HW}.
\renewcommand{\themainthm}{\ref{F-inversion}}
\begin{mainthm}
Let $(R, \m)$ be an F-finite regular local ring of characteristic $p>0$ and $I \subsetneq R$ an unmixed reduced ideal. Let $\a_1, \dots, \a_k$ be ideals of $R$ and $t_1, \dots, t_k > 0$ real numbers. We denote $S=R/I$.
If the pair $(S, (\a_1S)^{t_1}\cdots (\a_kS)^{t_k})$ is F-pure $($resp. strongly F-regular$)$, then the pair $(R, I\a_1^{t_1}\cdots \a_k^{t_k})$ is F-pure $($resp. purely F-regular$)$. 
\end{mainthm}

A pair $(R, \a^t)$ of a $\Q$-Gorenstein normal local ring $R$ of characteristic zero and an ideal $\a \subset R$ with real exponent $t > 0$ is of strongly F-regular type if and only if the pair $(\Spec R, t \cdot V(\a))$ is klt (cf. \cite[Theorem 6.8]{HY}), and if the pair $(R, \a^t)$ is of dense F-pure type (resp. purely F-regular type), then $(\Spec R, t \cdot V(\a))$ is lc (resp. plt) (Proposition \ref{plt}).
As an application of Theorem \ref{F-inversion}, via the above correspondence between singularities of pairs and ``F-singularities of pairs,'' we obtain a sort of Inversion of Adjunction of arbitrary codimension for klt pairs.

\renewcommand{\themainthm}{\ref{klt}}
\begin{mainthm}
Let $X$ be a non-singular variety over a field of characteristic zero and $Y=\sum_{i=1}^k t_iY_i$ a formal combination where $t_i > 0$ are real numbers and $Y_i \subsetneq X$ are closed subschemes. Let $Z \subsetneq X$ be a normal $\Q$-Gorenstein closed subvariety such that $Z \not\subset \cup_{i=1}^{k}Y_i$.
If the pair $(Z, Y|_Z)$ is klt, then the pair $(X, Y+Z)$ is plt near $Z$.
\end{mainthm}

On the other hand, the equivalence of lc pairs and pairs of dense F-pure type is still open (we do not know whether lc pairs are of dense F-pure type), hence we cannot derive the lc case from Theorem \ref{F-inversion} directly. Thus we need other techniques.

Hara \cite{Ha3} and Smith \cite{Sm2} independently proved that in a normal $\Q$-Gorenstein ring $R$ of characteristic zero, the multiplier ideal $\J(R)$ associated to the unit ideal corresponds to the test ideal $\tau(R)$ which plays a central role in the theory of tight closure.
Recently Hara and Yoshida \cite{HY} introduced a generalization of the test ideal, which is the ideal $\tau(\a^t)$ associated to a given ideal $\a$ with real exponent $t> 0$, and they extended Hara and Smith's result: the multiplier ideal $\J(\a^t)$ associated to a given ideal $\a$ with real exponent $t> 0$ in a normal $\Q$-Gorenstein ring $R$ of characteristic zero coincides, after reduction to characteristic $p \gg 0$, with the ideal $\tau(\a^t)$.

Since lc pairs are characterized by multiplier ideals (Lemma \ref{lc}), using some results about the ideal $\tau(\a^t)$ which are shown by a similar argument to the proof of Theorem \ref{F-inversion} (Corollary \ref{restriction}, etc.), we prove the lc case of Theorem \ref{klt} via the correspondence of $\J(\a^t)$ and $\tau(\a^t)$. 
\renewcommand{\themainthm}{\ref{inversion}}
\begin{mainthm}
Let $(X,Y)$ be a pair and $Z \subsetneq X$ a closed subvariety as in Theorem 4.1. 
If the pair $(Z, Y|_Z)$ is lc, then the pair $(X,Y+Z)$ is also lc near $Z$.
\end{mainthm}

The case where $Z$ is a divisor was proved by Ein, Musta{\c{t}}{\v{a}} and Yasuda \cite[Corollary 1.8]{EMY}, and we give another proof of their result by characteristic $p$ methods.

\section{Singularities of pairs and multiplier ideal sheaves}
We use the theory of singularities of pairs and multiplier ideal sheaves. Our main references are \cite{Ko} for the theory of singularities of pairs and \cite{La} for the theory of multiplier ideal sheaves. 
However, as we work with pairs of arbitrary codimension, we review some extensions to this setting of the definitions we need.

Let $X$ be a $\Q$-Gorenstein normal variety over a field of characteristic zero and $Y=\sum_{i=1}^k t_iY_i$ a formal combination where $t_i>0$ are real numbers and $Y_i \subsetneq X$ are closed subschemes defined by ideal sheaves $\a_i \subset \O_X$.
Suppose $f:\widetilde{X} \to X$ is a birational morphism from a normal variety $\widetilde{X}$ such that all ideal sheaves $\a_i\O_{\widetilde{X}}=\O_{\widetilde{X}}(-G_i)$ are invertible.
Let $K_X$ and $K_{\widetilde{X}}$ denote canonical divisors of $X$ and $\widetilde{X}$, respectively.
The relative canonical divisor $K_{\widetilde{X}/X}=K_{\widetilde{X}}-f^*K_X$ is a $\Q$-divisor supported on the exceptional locus $\mathrm{Exc}(f)$ of $f$. 
Then, for all the irreducible (not necessarily exceptional) divisors $E$ on $\widetilde{X}$, there are real numbers $a(E,X,Y)$ such that
$$K_{\widetilde{X}/X}-\sum_{i=1}^k t_iG_i = \sum_{E\textup{:arbitrary}} a(E,X,Y)E. $$
The number $a(E,X,Y)$ is called the discrepancy of $E$ with respect to the pair $(X,Y)$.

\begin{defn}\label{pair}
Let $(X,Y)$ be the same as above. 
\renewcommand{\labelenumi}{(\roman{enumi})}
A \textit{divisor over $(X,Y)$} means an irreducible divisor on some normal variety $\widetilde{X}$ with a birational morphism $f:\widetilde{X} \to X$ such that all $\a_i\O_{\widetilde{X}}$ are invertible.
\begin{enumerate}
\item We say that the pair $(X,Y)$ is \textit{log canonical} (or \textit{lc} for short) if $a(E,X,Y) \ge -1$ for all divisors $E$ over $(X,Y)$ (that is, $E$ runs through all the irreducible divisors of all birational morphisms $\widetilde{X} \to X$ such that $\a_i\O_{\widetilde{X}}$ is invertible for every $i=1, \dots, k$). 
\item We say that the pair $(X,Y)$ is \textit{Kawamata log terminal} (or \textit{klt} for short) if $a(E,X,Y) >-1$ for all divisors $E$ over $(X,Y)$. 
\item We say that the pair $(X,Y)$ is \textit{purely log terminal} (or \textit{plt} for short) if $a(E,X,Y) \ge -1$ for all divisors $E$ over $(X,Y)$ which dominate some irreducible component of $\cup_{i=1}^k Y_i$ and if $a(E,X,Y) > -1$ for the other divisors $E$ over $(X,Y)$. 
\item 
Suppose $f:\widetilde{X} \to X$ is a log resolution of the pair $(X, Y)$, that is, $f$ is a proper birational morphism with $\widetilde{X}$ non-singular such that all ideal sheaves $\a_i\O_{\widetilde{X}}=\O_{\widetilde{X}}(-G_i)$ are invertible, and $\cup_{i=1}^k \Supp G_i \cup \mathrm{Exc}(f)$ is a simple normal crossing divisor (Hironaka \cite{Hi} proved that log resolutions always exist).
Then the \textit{multiplier ideal sheaf} $\mathcal{J}(X, \a_1^{t_1} \cdots \a_k^{t_k})$ associated to $\a_1, \dots, \a_k$ with exponents $t_1, \dots, t_k$ is defined to be the ideal sheaf
$$\mathcal{J}(X, \a_1^{t_1} \cdots \a_k^{t_k})=f_*\mathcal{O}_{\widetilde{X}}(\lceil K_{\widetilde{X}/X}-\sum_{i=1}^k t_iG_i \rceil) \subseteq \mathcal{O}_X.$$
\end{enumerate}
\end{defn}

\begin{rem}\label{lcrem}
$(1)$ 
The Multiplier ideal sheaf does not depend on the choice of a log resolution $f:\widetilde{X} \to X$. 
The conditions (i) and (ii) of Definition \ref{pair} are checked by examining discrepancies of irreducible divisors on some log resolution $\widetilde{X}$ of the pair $(X,Y)$. 
The condition (iii) is checked by examining discrepancies of irreducible divisors on a log resolution $\widetilde{X}$ of $(X,Y)$ such that $\sum_{j=1}^mE_j$ is smooth, where $E_1, \dots, E_m$ are all irreducible divisors on $\widetilde{X}$ which dominate some irreducible component of $\cup_{i=1}^k Y_i$. 

$(2)$
The notions of klt, plt and lc pairs are local.
We say that the pair $(X,Y)$ is klt (resp. plt, lc) at a point $x \in X$ if the pair $(U, Y|_U)$ is klt (resp. plt, lc) for some neighborhood $U$ of $x$. 
Then the pair $(X,Y)$ is klt (resp. plt, lc) if and only if $(X,Y)$ is klt (resp. plt, lc) at every point $x \in X$.

$(3)$
The pair $(X, Y)$ is klt if and only if $\mathcal{J}(X, \a_1^{t_1} \cdots \a_k^{t_k})=\O_X$. 

$(4)$
If the pair $(X,Y)$ is lc, then $\mathcal{J}(X, \a_1^{t_1} \cdots \a_k^{t_k})$ is a radical ideal. 
\end{rem}

Like klt pairs, we can characterize lc pairs in terms of multiplier ideals.
\begin{lem}\label{lc}
Let $X$ be a $\Q$-Gorenstein normal variety over a field of characteristic zero and $Y=\sum_{i=1}^k t_iY_i$ a formal combination where $t_i > 0$ are real numbers and $Y_i \subsetneq X$ are closed subschemes defined by ideal sheaves $\a_i \subset \O_X$.
\begin{enumerate}
\item 
If there exists a nonzero ideal sheaf $\I \subseteq \O_X$ such that for every $0<s<1$, we have $\I \subset \J(X,\I^s\a_1^{st_1}\cdots \a_k^{st_k})$,
then the pair $(X,Y)$ is lc.
\item
Let $Z \subsetneq X$ be the non-klt locus of $X$ with the reduced induced closed subscheme structure and we denote by $\I_Z \subseteq \O_X$ the defining ideal of $Z$.
Then the pair $(X,Y)$ is lc if and only if $\I_Z \subset \J(X,\I_Z^s\a_1^{st_1}\cdots \a_k^{st_k})$ for all $0<s<1$.
\end{enumerate}
\end{lem}
\begin{proof}
$(1)$
Take a common log resolution $f: \widetilde{X} \to X$ of $\a_1, \dots, \a_k$ and $\I$ so that $\a_1\O_{\widetilde{X}}=\O_{\widetilde{X}}(-G_1)$, \dots, $\a_k\O_{\widetilde{X}}=\O_{\widetilde{X}}(-G_k)$ and $\I\O_{\widetilde{X}}=\O_{\widetilde{X}}(-F)$
 are invertible. The inclusion $\I \subset \J(X,\I^s\a_1^{st_1}\cdots \a_k^{st_k})$ implies that
$$ \lceil K_{\widetilde{X}/X}-sF-\sum_{i=1}^k st_iG_i \rceil \ge -F,$$
namely the coefficient of $K_{\widetilde{X}/X}+(1-s)F-s\sum_{i=1}^k t_iG_i$ in each irreducible component is greater than $-1$. Choosing $s$ very close to $1$, we see that the coefficient of $K_{\widetilde{X}/X}-\sum_{i=1}^k t_iG_i$ in each irreducible component is greater than or equal to $-1$. 

$(2)$
Assume that the pair $(X,Y)$ is lc.
Let $f:\widetilde{X} \to X$ be a common log resolution of $\a_1, \dots, \a_k$ and $\I_Z$ so that $\a_1\O_{\widetilde{X}}=\O_{\widetilde{X}}(-G_1)$, \dots, $\a_k\O_{\widetilde{X}}=\O_{\widetilde{X}}(-G_k)$ and $\I_Z\O_{\widetilde{X}}=\O_{\widetilde{X}}(-F_Z)$ are invertible.
Here we decompose $\lceil K_{\widetilde{X}/X}\rceil=K_{+}-K_{-}$, where $K_+$ and $K_-$ are effective divisors on $\widetilde{X}$ which have no common component.
Notice that $K_-$ is reduced or zero, because the pair $(X, Y)$ is lc and in particular $X$ has only log canonical singularities.
By the definition of $Z$ and Remark \ref{lcrem} $(4)$, $\I_Z$ is equal to the multiplier ideal sheaf $\J(X, \O_X)=f_*\O_{\widetilde{X}}(\lceil K_{\widetilde{X}/X} \rceil)=f_*\O_{\widetilde{X}}(-K_-)$, hence $F_Z \ge K_-$.
Since the pair $(X,Y)$ is lc, that is, the coefficient of $K_{\widetilde{X}/X}-\sum_{i=1}^k t_iG_i$ in each irreducible component is  not less than $-1$,
$$\lceil K_{\widetilde{X}/X}+(1-s)F_Z-s\sum_{i=1}^k t_iG_i \rceil \ge 0$$ for every $0<s<1$.
This implies that $\I_Z \subset \J(X,\I_Z^s\a_1^{st_1}\cdots \a_k^{st_k})$ for all $0<s<1$.
\end{proof}

Nakayama \cite{Na} introduced the notion of admissible pairs as an analog of klt pairs in the absence of $\Q$-Gorensteinness.
\begin{defn}
{\rm (\cite[Definition A.1.2]{Na})}
Let $(X, \Delta)$ be a pair of a normal variety $X$ over a field of characteristic zero and an effective $\R$-divisor $\Delta$ with $\lfloor \Delta \rfloor=0$. The pair $(X, \Delta)$ is called {\it strictly admissible} if there exist a resolution of singularities $f:\widetilde{X} \to X$ and a $\Q$-divisor $E$ on $\widetilde{X}$ satisfying the following conditions:
\renewcommand{\labelenumi}{(\roman{enumi})}
\begin{enumerate}
\item $\Supp E$ is a normal crossing divisor.
\item $\lceil E \rceil$ is an effective $f$-exceptional divisor.
\item $-f_*E \ge \Delta$.
\item $E-K_{\widetilde{X}}$ is $f$-ample.
\end{enumerate}
The pair $(X, \Delta)$ is called {\it admissible} if there is an open covering $\{U_{\lambda}\}$ of $X$ such that $(U_{\lambda}, \Delta|_{U_{\lambda}})$ is strictly admissible for any $\lambda$.
\end{defn}

\begin{rem}
The notion of admissible pairs is a generalization of the notion of klt pairs: under the condition that $K_X+\Delta$ is $\R$-Cartier, $(X,\Delta)$ is admissible if and only if $(X, \Delta)$ is klt (see \cite[Lemma A.1.7]{Na}). Refer to \cite{Na} for other properties of admissible pairs.
\end{rem}

Hara and Watanabe \cite{HW} defined the notion of strong F-regularity for a pair $(R,D)$ of a normal ring $R$ and an effective $\R$-divisor $D$ on $\Spec R$ (refer to \cite{HW} and \cite{Ta} for the definition and basic properties, and see also Remark \ref{divisor}).
We say that a pair $(X, \Delta)$ of a normal variety $X$ over a field of characteristic zero and an effective $\R$-divisor $\Delta$ on $X$ is of {\it strongly F-regular type} if  the pair $(\O_{X,x}, \Delta_x)$ is of strongly F-regular type for all points $x \in X$.
By virtue of \cite[Corollary 3.4]{Ta}, admissible pairs are of strongly F-regular type. 
\begin{lem}\label{admissible}
Let $X$ be a normal variety over a field of characteristic zero and $\Delta$ an effective $\R$-divisor on $X$.
If the pair $(X, \Delta)$ is admissible, then $(X, \Delta)$ is of strongly F-regular type.
\end{lem}
\begin{proof}
Fix any point $x \in X$.
If the pair $(X, \Delta)$ is admissible, then by \cite[Lemma A.1.3]{Na}, there exist an open neighborhood $U$ of $x$ and an effective $\Q$-divisor $\Delta'$ on $U$ such that $\Delta' \ge \Delta|_U$ and $(U, \Delta')$ is klt. It follows from \cite[Corollary 3.4]{Ta} that the pair $(U, \Delta')$ is of strongly F-regular type. Since $\Delta' \ge \Delta|_U$, by \cite[Proposition 2.2]{HW}, the pair $(U, \Delta|_U)$ is also of strongly F-regular type. Since strong F-regularity is a local property, we conclude that $(X, \Delta)$ is of strongly F-regular type.
\end{proof}

\section{Brief review on a generalization of test ideals}
The notion of tight closure was introduced by Hochster and Huneke \cite{HH1} by using the Frobenius map in characteristic $p>0$.
The test ideal $\tau(R)$ of a ring $R$ of characteristic $p>0$ is an important object in the theory of tight closure.
Hara and Yoshida \cite{HY} introduced a generalization of the test ideal, which is the ideal $\tau(\a^t)$ associated to a given ideal $\a$ with real exponent $t > 0$, and investigated properties of that ideal.
In this section, we briefly review the definition and fundamental properties of the ideal $\tau(\a^t)$ which we need later.
Refer to \cite{HY} for the proofs.

In this paper, all rings are reduced commutative rings with unity. 
For a ring $R$, we denote by $R^{\circ}$ the set of elements of $R$ which are not in any minimal prime ideal. 
Let $R$ be a ring of characteristic $p >0$ and $F\colon R \to R$ the Frobenius map which sends $x \in R$ to $x^p \in R$. 
For an integer $e > 0$, the ring $R$ viewed as an $R$-module via the $e$-times iterated Frobenius map $F^e \colon R \to R$ is denoted by ${}^e\! R$. 
Since $R$ is assumed to be reduced, we can identify $F^e \colon R \to {}^e\! R$ with the natural inclusion map $R \hookrightarrow R^{1/p^e}$. We say that $R$ is {\it F-finite} if ${}^1
\! R$ (or $R^{1/p}$) is a finitely generated $R$-module. 
We always assume that all rings of characteristic $p>0$ are F-finite throughout this paper. 

Let $R$ be a ring of characteristic $p > 0$ and $M$ an $R$-module. 
For each integer $e > 0$, we denote $\F^e(M) = \F_R^e(M) := M \otimes_R {}^e\! R$ and 
regard it as an $R$-module by the action of $R$ on ${}^e\! R$ from the right. 
Then we have the induced $e$-times iterated Frobenius map $F^e \colon M \to 
\F^e(M)$. 
The image of $z \in M$ via this map is denoted by $z^q:= F^e(z) 
\in \F^e(M)$, where $q=p^e$. For an $R$-submodule $N$ of $M$, we denote by $N^{[q]}_M$ the 
image of the induced map $\F^e(N) \to \F^e(M)$. 
If $M=R$ and $N=I \subseteq R$, then we write $I^{[q]}=I^{[q]}_R$.

When $R$ is regular, the Frobenius functor $\F^e_R$ has a nice property.
\begin{lem}\label{flat}
{\rm (\cite{Ku})}
Let $R$ be a regular ring of characteristic $p>0$.
Then for every integer $e > 0$, $\F^e_R$ is an exact functor, that is, $R \hookrightarrow R^{1/p^e}$ is flat.
\end{lem}

Now we recall the definition of $\a_1^{t_1} \cdots \a_k^{t_k}$-tight closure.
\begin{defn}\label{ta}
{\rm (\cite[Definition 6.1]{HY})}
Let $\a_1, \dots, \a_k$ be ideals of a ring $R$ of characteristic $p>0$ and $N \subseteq M$ $R$-modules. 
Given real numbers $t_1, \dots, t_k > 0$, the {\it $\a_1^{t_1} \cdots \a_k^{t_k}$-tight closure} $N^{*\a_1^{t_1}\cdots \a_k^{t_k}}_M$ of $N$ in $M$ is defined to be the submodule of $M$ consisting of all elements $z \in M$ for which there exists $c \in R^{\circ}$ such that 
$$cz^q\a_1^{\lceil t_1q \rceil} \cdots \a_k^{\lceil t_kq \rceil} \subseteq N^{[q]}_M$$ 
for all large $q = p^e$, where $\lceil t_iq \rceil$ is the least integer which is greater than or equal to $t_iq$ for every $i=1, \dots, k$. 
If $N=N^{*\a_1^{t_1}\cdots \a_k^{t_k}}_M$, then we say that $N$ is $\a_1^{t_1} \cdots \a_k^{t_k}$-{\it tightly closed} in $M$.
The $\a_1^{t_1} \cdots \a_k^{t_k}$-tight closure of an ideal $I \subseteq R$ is just defined by $I^{*\a_1^{t_1} \cdots \a_k^{t_k}} = I^{*\a_1^{t_1} \cdots \a_k^{t_k}}_R$.
\end{defn}

\begin{rem}
When $\a$ is the unit ideal, $\a^t$-tight closure is the same as classical tight closure. Refer to Hochster and Huneke's original paper \cite{HH1} for the classical tight closure theory.  
\end{rem}

In order to handle the $\a_1^{t_1} \cdots \a_k^{t_k}$-tight closure operation, the notion of $\a_1^{t_1} \cdots \a_k^{t_k}$-test elements is quite useful.
\begin{defn}\label{testdef}
{\rm (\cite[Definition 6.3]{HY})}
Let $\a_1, \dots, \a_k$ be ideals of a ring $R$ of characteristic $p > 0$ and $t_1, \dots, t_k > 0$ real numbers. An element $d \in R^{\circ}$ is called an {\it $\a_1^{t_1} \cdots \a_k^{t_k}$-test element} if for every finitely generated $R$-module 
$M$ and $z \in M$, the following holds: $z \in 0_M^{*\a_1^{t_1} \cdots \a_k^{t_k}}$ if and only 
if $dz^q\a_1^{\lceil t_1q \rceil} \cdots \a_k^{\lceil t_kq \rceil}= 0$ for all powers $q = p^e$ of $p$.
\end{defn}

An $\a_1^{t_1} \cdots \a_k^{t_k}$-test element exists in almost every ring. Particularly in case a ring $R$ is regular, the unit element $1$ is an $\a_1^{t_1} \cdots \a_k^{t_k}$-test element for any ideals $\a_1, \dots, \a_k \subseteq R$ and any real numbers $t_1, \dots, t_k > 0$.
\begin{thm}\label{exist}
{\rm (\cite[Theorem 6.4]{HY})}
Let $R$ be an F-finite reduced ring of characteristic $p>0$ and $c \in R^{\circ}$ an element such that the localized ring $R_c$ is strongly F-regular $($as for the definition of strongly F-regular rings, see Remark \ref{divisor}$)$. Then some power $c^n$ of c is an $\a_1^{t_1} \cdots \a_k^{t_k}$-test element for all ideals $\a_1, \dots, \a_k \subseteq R$ and all real numbers $t_1, \dots, t_k > 0$.
\end{thm}

Now using the $\a_1^{t_1} \dots \a_k^{t_k}$-tight closure of the zero submodule, we define the ideal $\tau(R, \a_1^{t_1}\cdots \a_k^{t_k})$.
The following is a generalization of \cite[Definition 8.22, Proposition 8.23 and 8.25]{HH1} and \cite{AM}. 
\begin{defthm}\label{taudef}
{\rm (\cite[Definition-Theorem 6.5]{HY})}
Let $R$ be an excellent reduced ring of characteristic $p > 0$, $\a_1, \dots, \a_k$ ideals of $R$ and $t_1, \dots, t_k > 0$ real numbers. 
Let $E =\bigoplus_{\m} E_R(R/\m)$ be the direct sum, taken over all maximal ideals $\m$ of $R$, of the injective hulls of the residue fields $R/\m$. 
Then the following ideals are equal to each other 
and we denote them by $\tau(R,\a_1^{t_1}\cdots \a_k^{t_k})$. 
\begin{enumerate}
\renewcommand{\labelenumi}{\textup{(\roman{enumi})}}
\item $\displaystyle\bigcap_M \Ann_R(0^{*\a_1^{t_1}\cdots \a_k^{t_k}}_M)$, where $M$ runs through all finitely generated $R$-modules. 

\item $\displaystyle\bigcap_{M\subseteq E} \Ann_R(0^{*\a_1^{t_1}\cdots \a_k^{t_k}}_M)$, where $M$ runs through all finitely generated submodules of $E$. 
\item $\displaystyle\bigcap_{J\subseteq R} (J:J^{*\a_1^{t_1}\cdots \a_k^{t_k}})$, where $J$ runs through all ideals of $R$. 
\end{enumerate}
Moreover if $R$ is a normal $\Q$-Gorenstein local ring, then the above three ideals also coincide with the following ideal.
\renewcommand{\labelenumi}{\textup{(\roman{enumi})}}
\begin{enumerate}
\setcounter{enumi}{3}
\item $\displaystyle\Ann_R(0^{*\a_1^{t_1}\cdots \a_k^{t_k}}_E)$.
\end{enumerate}
\end{defthm}

\begin{rem}\label{testrem} 
In the case where $\a = R$ is the unit ideal, the ideal $\tau(\a) 
= \tau(R)$ is called the test ideal of $R$ and an $R$-test element 
is nothing but a test element as defined in \cite{HH1}. 
In this case, $\tau(R) \cap R^{\circ}$ is exactly equal to the set of test elements. 
However, $\tau(\a^t) \cap R^{\circ}$ is not equal to the set of $\a^t$-test elements in general. 
We also remark that by \cite[Corollary 2.4]{HT}, a test element is an $\a^t$-test element for all ideals $\a \subseteq R$ and all real numbers $t>0$ under some mild condition. 
\end{rem}

Let $R$ be an algebra essentially of finite type over a field $k$ of characteristic zero, $\a_1, \dots, \a_k \subset R$ ideals and $t_1, \dots, t_k>0$ real numbers. One can choose a finitely generated $\Z$-subalgebra $A$ of $k$ and a subalgebra $R_A$ of $R$ essentially of finite type over $A$ such that the natural map $R_A \otimes_A k \to R$ is an isomorphism and ${\a_A}_iR=\a_i$ where ${\a_A}_i=\a_i \cap R_A \subset R_A$ for all $1 \le i \le k$. 
Given a closed point $s \in \Spec A$ with residue field $\kappa=\kappa(s)$, we denote the corresponding fibers over $s$ by $R_{\kappa}, {\a_{\kappa}}_1, \dots, {\a_{\kappa}}_k$ (cf. Definition \ref{char 0}).
Then we refer to such $(\kappa, R_{\kappa}, {\a_{\kappa}}_1, \dots, {\a_{\kappa}}_k)$ for a general closed point $s \in \Spec A$ with residue field $\kappa=\kappa(s)$ of sufficiently large characteristic $p \gg 0$ as ``{\it reduction to characteristic $p \gg 0$}'' of $(k,R, \a_1, \dots, \a_k)$, and the pair $(R_{\kappa}, {\a_{\kappa}}_1^{t_1} \cdots {\a_{\kappa}}_k^{t_k})$ inherits the properties possessed by the original pair $(R, \a_1^{t_1} \cdots \a_k^{t_k})$.
Furthermore, given a log resolution $f:\widetilde{X} \to X=\Spec R$ of $(X, \a_1^{t_1} \dots \a_k^{t_k})$, we can reduce this entire setup to characteristic $p \gg 0$. 

Hara and Yoshida extend the result of Hara \cite{Ha2} and Smith \cite{Sm2} to that for pairs. 
\begin{thm}\label{multiplier}
{\rm (\cite[Theorem 6.8]{HY})}
Let $(R,\m)$ be a $\Q$-Gorenstein normal local ring essentially of finite type over a field of characteristic zero. 
Let $\a_1, \dots, \a_k \subseteq R$ be nonzero ideals and $t_1, \dots, t_k > 0$ real numbers. 
Then, after reduction to characteristic $ p \gg 0$, 
$$\tau(R, \a_1^{t_1} \cdots \a_k^{t_k})=\J(\Spec R, \a_1^{t_1} \cdots \a_k^{t_k}).$$
\end{thm}
By virtue of the above theorem, we can study local properties of the multiplier ideal $\J(\Spec R, \a_1^{t_1} \cdots \a_k^{t_k})$ by examining those of the ideal $\tau(R, \a_1^{t_1}\cdots \a_k^{t_k})$. By using this, we prove a sort of Inversion of Adjunction of arbitrary codimension for lc pairs in section $4$.

\section{F-singularities of pairs}
Hara and Watanabe \cite{HW} defined the notions of F-regularity and F-purity for a pair $(R, \Delta)$ of a normal ring $R$ of characteristic $p>0$ and an effective $\R$-divisor $\Delta$ on $\Spec R$. 
In this section, we introduce the notions of F-regularity and F-purity for a pair $(R, \a^t)$ of a ring $R$ of characteristic $p>0$ and an ideal $\a \subset R$ with real exponent $t > 0$, and examine these properties.

\begin{defn}\label{ideal}
Let $\a_1, \dots, \a_k$ be ideals of an F-finite reduced ring $R$ of characteristic $p>0$ and $t_1, \dots, t_k > 0$ real numbers.
\renewcommand{\labelenumi}{(\roman{enumi})}
\begin{enumerate}
\item The pair $(R, \a_1^{t_1} \cdots \a_k^{t_k})$ is said to be {\it F-pure} (or $R$ is said to be {\it F-pure} with respect to $\a_1, \dots, \a_k$ and $t_1, \dots, t_k$) if for all large $q=p^e$, there exists an element $d \in \a_1^{\lfloor t_1(q-1) \rfloor} \cdots \a_k^{\lfloor t_k(q-1) \rfloor}$ such that $d^{1/q}R \hookrightarrow R^{1/q}$ splits as an $R$-module homomorphism.
\item The pair $(R, \a_1^{t_1} \cdots \a_k^{t_k})$ is said to be {\it strongly F-regular} (or $R$ is said to be {\it strongly F-regular} with respect to $\a_1, \dots, \a_k$ and $t_1, \dots, t_k$) if for every $c \in R^{\circ}$, there exist $q=p^e$ and $d \in \a_1^{\lfloor t_1(q-1) \rfloor} \cdots \a_k^{\lfloor t_k(q-1) \rfloor}$ such that $(cd)^{1/q}R \hookrightarrow R^{1/q}$ splits as an $R$-module homomorphism. 
\item The pair $(R, \a_1^{t_1} \cdots \a_k^{t_k})$ is said to be {\it purely F-regular} (or $R$ is said to be {\it purely F-regular} with respect to $\a_1, \dots, \a_k$ and $t_1, \dots, t_k$) if for every element $c \in R^{\circ}$ which is not in any minimal prime ideal of $\a_1 \cdots \a_k$, there exists $q'$ such that for all $q=p^e \ge q'$ and for some $d \in \a_1^{\lfloor t_1(q-1) \rfloor} \cdots \a_k^{\lfloor t_k(q-1) \rfloor}$, $(cd)^{1/q}R \hookrightarrow R^{1/q}$ splits as an $R$-module homomorphism. 
\end{enumerate}
\end{defn}

\begin{rem}\label{divisor}
$(1)$
By definition, the pure F-regularity with respect to the unit ideal is equivalent to the strong F-regularity with respect to the unit ideal.
A ring $R$ is F-pure (resp. strongly F-regular) if and only if $R$ is F-pure (resp. purely F-regular) with respect to the unit ideal $R$. 
Refer to \cite{HR}, \cite{HH1} and \cite{HH2} for properties of F-pure and strongly F-regular rings.

$(2)$
We have analogous notions of strong F-regularity and F-purity for a pair $(R,\Delta)$ of a normal ring $R$ of characteristic $p>0$ and an effective $\R$-divisor $\Delta$ on $\Spec R$. See \cite{HW} and \cite{Ta} for details. If $\a_i=x_iR$ and $\Delta=\sum_{i=1}^k t_i \cdot \Div(x_i)$ with $x_i \in R^{\circ}$ and $t_i \in \R_{> 0}$ for $1 \le i \le n$, then the F-purity (resp. strong F-regularity) of the pair $(R, \a_1^{t_1} \cdots \a_k^{t_k})$ is the same as the F-purity (resp. strong F-regularity) of the pair $(R,\Delta)$. 

$(3)$
We can generalize the notion of global F-regularity introduced in \cite{HWY}, \cite{Sm3} and consider the global version of strong F-regularity with respect to ideal sheaves $\I_1, \dots, \I_k$ and real numbers $t_1, \dots, t_k>0$, by using the $e$-Frobenius splitting along an effective Cartier divisor $D$ such that $f^*D \ge \sum_{i=1}^k \lfloor t_i(q-1) \rfloor F_i$ where $f:Y \to X$ is the normalized blowing-up along $\I_1, \dots, \I_k$ so that $\I_i\O_Y=\O_Y(-F_i)$ are invertible. 
Refer to \cite{Sm3} for Frobenius splitting along a divisor and global F-regularity. 
\end{rem}

We collect some basic properties of ``F-singularities of pairs.''
\begin{prop}\label{basic}
Let $\a_1, \dots, \a_k$ be ideals of an F-finite reduced ring $R$ of characteristic $p>0$ and $t_1, \dots, t_k > 0$ real numbers.
\begin{enumerate}
\item If the pair $(R, \a_1^{t_1} \cdots \a_k^{t_k})$ is strongly F-regular $($resp. F-pure$)$, then so is the pair $(R, \b_1^{s_1} \cdots \b_k^{s_k})$ for every ideal $\b_i \supseteq \a_i$ and every real number $0 <s_i \le t_i$ for all $1 \le i \le k$.
\item A strongly F-regular pair is purely F-regular, and a purely F-regular pair is F-pure. 
\item The pair $(R,\a_1^{t_1} \cdots \a_k^{t_k})$ is strongly F-regular if and only if for every $c \in R^{\circ}$, there exists $q'$ such that for all $q=p^e \ge q'$ and for some $d \in \a_1^{\lceil t_1q \rceil} \cdots \a_k^{\lceil t_kq \rceil}$, $(cd)^{1/q}R \hookrightarrow R^{1/q}$ splits as an $R$-module homomorphism.
\end{enumerate}
\end{prop}
\begin{proof}
$(1)$ is obvious. Since $(2)$ immediately follows from $(3)$, we prove only $(3)$. The sufficiency is clear. To show the necessity, choose an element $a \in R^{\circ}$ such that $a\a_1^{\lfloor t_1(q-1) \rfloor} \cdots \a_k^{\lfloor t_k(q-1) \rfloor} \subset \a_1^{\lceil t_1q \rceil} \cdots \a_k^{\lceil t_kq \rceil}$ for every $q=p^e$. 
Since the pair $(R,\a_1^{t_1} \cdots \a_k^{t_k})$ is strongly F-regular, for every element $c \in R^{\circ}$, there exist a power $q'$ of $p$ and an element $d' \in \a_1^{\lfloor t_1(q'-1) \rfloor} \cdots \a_k^{\lfloor t_k(q'-1) \rfloor}$ such that $R \xrightarrow{(acd')^{1/q'}} R^{1/q'}$ splits as an $R$-module homomorphism. 
Since $R$ is in particular F-pure by $(1)$, the map $R^{1/q'} \hookrightarrow R^{1/qq'}$ splits for all $q=p^e$.
Therefore the composite map $R \xrightarrow{(acd')^{1/q'}} R^{1/q'} \hookrightarrow R^{1/qq'}$ also splits for all $q=p^e$.
This map is factorized into $$R \xrightarrow{c^{1/qq'}(ad')^{q/qq'}} R^{1/qq'} \xrightarrow{c^{(q-1)/qq'}} R^{1/qq'}$$ and $(ad')^q \in (\a_1^{\lceil t_1q' \rceil} \cdots \a_k^{\lceil t_kq' \rceil})^q \subset \a_1^{\lceil t_1qq' \rceil} \cdots \a_k^{\lceil t_kqq' \rceil}$, hence the map $R \xrightarrow{(cd)^{1/q}} R^{1/q}$ splits for all  powers $q=p^e \ge q'$ and for some element $d \in \a_1^{\lceil t_1q \rceil} \cdots \a_k^{\lceil t_kq \rceil}$.
\end{proof}

\begin{lem}\label{injective}
Let $(R,\m)$ be an F-finite reduced local ring of characteristic $p>0$. Let $\a_1, \dots, \a_k$ be ideals of $R$ and $t_1, \dots, t_k > 0$ real numbers. We denote by $E_R$ the injective hull of the residue field $R/\m$ and by $F^e:E_R \to E_R \otimes_R R^{1/q}$ the induced $e$-times iterated Frobenius map on $E_R$.
\begin{enumerate}
\item
The pair $(R,\a_1^{t_1} \cdots \a_k^{t_k})$ is F-pure if and only if for all large $q=p^e$, there exists an element $d \in \a_1^{\lfloor t_1(q-1) \rfloor} \cdots \a_k^{\lfloor t_k(q-1) \rfloor}$ such that $dF^e:E_R \to E_R \otimes_R R^{1/q}$ is injective.
\item 
The pair $(R, \a_1^{t_1} \cdots \a_k^{t_k})$ is purely F-regular if and only if for every element $c \in R^{\circ}$ which is not in any minimal prime ideal of $\a_1 \cdots \a_k$, there exists $q'$ such that for all $q=p^e \ge q'$ and for some $d \in \a_1^{\lfloor t_1(q-1) \rfloor} \cdots \a_k^{\lfloor t_k(q-1) \rfloor}$, $cdF^e:E_R \to E_R \otimes_R R^{1/q}$ is injective.
\item
The pair $(R, \a_1^{t_1} \cdots \a_k^{t_k})$ is strongly F-regular if and only if for every element $c \in R^{\circ}$, there exist $q=p^e$ and $d \in \a_1^{\lfloor t_1(q-1) \rfloor} \cdots \a_k^{\lfloor t_k(q-1) \rfloor}$ such that $cdF^e:E_R \to E_R \otimes_R R^{1/q}$ is injective. 
\end{enumerate}
\end{lem}
\begin{proof}
This follows from a standard argument about the Matlis duality (cf. \cite[Proposition 2.4]{HW}).
Fix an element $d \in R^{\circ}$.
The map $d^{1/q}R \hookrightarrow R^{1/q}$ splits as an $R$-linear map if and only if the map
$$\Hom_R(R^{1/q}, R) \xrightarrow{d^{1/q}} \Hom_R(R^{1/q}, R) \to \Hom_R(R,R)=R$$
is surjective. 
By the local duality, the Matlis dual of $\Hom_R(R^{1/q}, R)$ is isomorphic to the $R$-module $E_R \otimes_R R^{1/q}$.
Therefore the surjectivity of the above map is equivalent to the injectivity of the map
$$dF^e:E_R \to E_R \otimes_R R^{1/q} \xrightarrow{d^{1/q}} E_R \otimes_R R^{1/q}.$$
Thus the assertion follows.
\end{proof}

\begin{cor}\label{F-regular}
Under the assumption of Lemma \ref{injective}, $(R, \a_1^{t_1} \cdots \a_k^{t_k})$ is a strongly F-regular pair if and only if $0_{E_R}^{*\a_1^{t_1} \cdots \a_k^{t_k}}=0$. In particular when $R$ is a normal $\Q$-Gorenstein local ring, this is equivalent to the condition that $\tau(R, \a_1^{t_1} \cdots \a_k^{t_k})=R$.
\end{cor}
\begin{proof}
Assume that the pair $(R,\a_1^{t_1} \cdots \a_k^{t_k})$ is strongly F-regular and fix any element $z \in 0_{E_R}^{*\a_1^{t_1} \cdots \a_k^{t_k}}$. By definition, there exists $c \in R^{\circ}$ such that $cdF^e(z)=0$ for all $q=p^e \gg 0$ and for all $d \in \a_1^{\lceil t_1q \rceil} \cdots \a_k^{\lceil t_kq \rceil}$. 
Since the pair $(R,\a_1^{t_1} \cdots \a_k^{t_k})$ is strongly F-regular, by Proposition \ref{basic} $(3)$ and Lemma \ref{injective}, there exist $q=p^e \gg 0$ and $d \in \a_1^{\lceil t_1q \rceil} \cdots \a_k^{\lceil t_kq \rceil}$ such that $cdF^e:E_R \to E_R \otimes_R R^{1/q}$ is injective, whence $z=0$.

Conversely suppose that $0_{E_R}^{*\a_1^{t_1} \cdots \a_k^{t_k}}=0$ and fix any element $c \in R^{\circ}$. 
If $z$ is a generator of the socle $(0:\m)_{E_R}$ of $E_R$, then by assumption, there exist $q=p^e$ and $d \in  \a_1^{\lceil t_1q \rceil} \cdots \a_k^{\lceil t_kq \rceil}$ such that $cdF^e(z) \ne 0$, that is, $cdF^e$ is injective on $(0:\m)_{E_R}$.
Since $E_R$ is an essential extension of $(0:\m)_{E_R}$, $cdF^e$ itself is injective. By Lemma \ref{injective}, this implies that the pair $(R,\a_1^{t_1} \cdots \a_k^{t_k})$ is strongly F-regular.
\end{proof}

\begin{lem}\label{test}
Let $R$ be an F-finite reduced ring of characteristic $p>0$ and $c \in R^{\circ}$ an element such that the localization $R_c$ with respect to $c$ is strongly F-regular. Then, for any ideals $\a_1, \dots, \a_k \subseteq R$ and for any real numbers $t_1, \dots, t_k>0$, the pair $(R, \a_1^{t_1}\cdots \a_k^{t_k})$ is strongly F-regular if and only if there exist $q=p^e$ and $d \in \a_1^{\lceil t_1q \rceil} \cdots \a_k^{\lceil t_kq \rceil}$ such that $(cd)^{1/q}R \hookrightarrow R^{1/q}$ splits as an $R$-linear map.
\end{lem}
\begin{proof}
By Proposition \ref{basic} $(3)$, the necessity is obvious. Thus we conversely suppose that there exist $q=p^e$ and $d \in \a_1^{\lceil t_1q \rceil} \cdots \a_k^{\lceil t_kq \rceil}$ such that $(cd)^{1/q}R \hookrightarrow R^{1/q}$ splits as an $R$-linear map, and fix any element $c' \in R^{\circ}$. By \cite[Remark 3.2]{HH2}, for some $q'=p^{e'}$ and $q''=p^{e''}$, there exists an $R$-module homomorphism
$$R^{1/q'} \to R, \quad {c'}^{1/q'} \mapsto c^{q''}.$$
Taking $qq''$-th roots, we obtain an $R^{1/qq''}$-linear map
$$R^{1/qq'q''} \to R^{1/qq''}, \quad {c'}^{1/qq'q''}d^{1/q} \mapsto (cd)^{1/q}.$$
Since by assumption it is easily shown that $R$ is F-pure, there exists an $R^{1/q}$-linear map
$$R^{1/qq''} \to R^{1/q}, \quad (cd)^{1/q} \mapsto (cd)^{1/q}.$$
By composing these maps, we have the following $R$-module homomorphism 
\begin{align*}
R^{1/qq'q''} &\longrightarrow R^{1/qq''} \longrightarrow R^{1/q} \longrightarrow R.\\ 
{c'}^{1/qq'q''}d^{q'q''/qq'q''} &\longmapsto (cd)^{1/q} \longmapsto (cd)^{1/q} \mapsto 1.
\end{align*} 
Since $d^{q'q''} \in \a_1^{\lceil t_1qq'q'' \rceil} \cdots \a_k^{\lceil t_kqq'q'' \rceil}$, the pair $(R, \a_1^{t_1} \cdots \a_k^{t_k})$ is strongly F-regular. 
\end{proof}

The notions of F-regularity and F-purity are also defined for a pair $(R, \a^t)$ of a ring $R$ of characteristic zero and an ideal $\a \subset R$ with real exponent $t > 0$.
\begin{defn}\label{char 0}
Let $R$ be a reduced algebra essentially of finite type over a field $k$ of characteristic zero, $t_1, \dots, t_k > 0$ real numbers, and $\a_1, \cdots, \a_k$ ideals of $R$.
The pair $(R, \a_1^{t_1} \cdots \a_k^{t_k})$ is said to be of {\it F-pure type} (resp. {\it strongly F-regular type}, {\it purely F-regular type}) if there exist a finitely generated $\Z$-subalgebra $A$ of $k$ and a reduced subalgebra $R_A$ of $R$ essentially of finite type over $A$ which satisfy the following conditions:
\renewcommand{\labelenumi}{(\roman{enumi})}
\begin{enumerate}
\item $R_A$ is flat over $A$, $R_A \otimes_A k = R$ and ${\a_A}_iR= \a_i$ where ${\a_A}_i=\a_i \cap R_A \subset R_A$ for all $i=1, \dots, k$.
\item The pair $(R_{\kappa}, {\a_{\kappa}}_1^{t_1} \cdots {\a_{\kappa}}_k^{t_k})$ is F-pure (resp. strongly F-regular, purely F-regular) for every closed point $s$ in a dense open subset of $\Spec A$, where $\kappa=\kappa(s)$ denotes the residue field of $s \in \Spec A$, $R_{\kappa}=R_A \otimes_A \kappa(s)$ and ${\a_{\kappa}}_i={\a_A}_i R_{\kappa} \subset R_{\kappa}$ for every $i=1, \dots, k$.
\end{enumerate}
The pair $(R, \a_1^{t_1} \cdots \a_k^{t_k})$ is said to be of {\it dense F-pure type} if in the above condition (ii) ``dense open" is replaced by ``dense." 
\end{defn}

By virtue of Theorem \ref{multiplier} and Corollary \ref{F-regular}, a pair $(R, \a^t)$ of a $\Q$-Gorenstein normal local ring $R$ of characteristic zero and an ideal $\a \subset R$ with real exponent $t > 0$ is of strongly F-regular type if and only if the pair $(\Spec R, t \cdot V(\a))$ is klt. We show that if the pair $(R, \a^t)$ is of purely F-regular type (resp. dense F-pure type), then $(\Spec R, t \cdot V(\a))$ is plt (resp. lc). 
\begin{prop}\label{plt}
Let $(R, \m)$ be a $\Q$-Gorenstein normal local ring essentially of finite type over a field of characteristic zero and write $X=\Spec R$. 
Let $Y=\sum_{i=1}^k t_iY_i$ be a formal combination where $t_i > 0$ are real numbers and $Y_i \subsetneq X$ are closed subschemes defined by nonzero ideals $\a_i \subset R$.
If the pair $(R, \a_1^{t_1} \cdots \a_k^{t_k})$ is of dense F-pure type $($resp. purely F-regular type, strongly F-regular type$)$, then the pair $(X, Y)$ is lc $($resp. plt, klt$)$.
\end{prop}
\begin{proof}
The proof is essentially the same as that in \cite[Theorem 3.3]{HW}. We only prove the purely F-regular case, because the other cases also follow from a similar argument.

Let $f:\widetilde{X} \to X$ be an arbitrary proper birational morphism with $\widetilde{X}$ normal such that all ideal sheaves $\a_i\O_{\widetilde{X}}=\O_{\widetilde{X}}(-G_i)$ are invertible. Then there are finitely many irreducible divisors $E_j$ on $\widetilde{X}$ such that
$$K_{\widetilde{X}} \underset{\text{$\Q$-lin.}}{\sim} f^*K_X + \sum_{j=1}^n a_jE_j+\sum_{i=1}^k t_iG_i, $$
where $a_j$ are real numbers chosen as $\sum_{j=1}^n a_jE_j+\sum_{i=1}^k t_iG_i$ is an $f$-exceptional divisor. 
Suppose that $E_1, \dots, E_m$ $(m \le n)$ are all divisors which dominate some irreducible component of $\cup_{i=1}^k Y_i$. 
Considering the reduction to characteristic $p \gg 0$, we may assume that $R$, $\a_1, \dots, \a_k$, $f$, etc. are defined over a field of characteristic $p>0$.

Assume that the pair $(R, \a_1^{t_1} \cdots \a_k^{t_k})$ is purely F-regular. Choose some element $c \in R^{\circ}$ which is not in any minimal prime ideal of $\a_1 \cdots \a_k$ such that $v_{E_j}(c)$ is not less than the coefficient of $\sum_{i=1}^k G_i$ in $E_j$ for every $j=m+1,\dots,n$, where $v_{E_j}$ is the valuation of $E_j$.
Then by definition, for sufficiently large $q=p^e$ and for some $d \in \a_1^{\lfloor t_1(q-1) \rfloor} \cdots \a_k^{\lfloor t_k(q-1) \rfloor}$, there exists an $R$-linear map $\psi:R^{1/q} \to R$ sending $(cd)^{1/q}$ to $1$. 
Let $\phi=\psi \circ (cd)^{1/q} \in \Hom_R(R^{1/q},R)$.
Via the isomorphism $\Hom_R(R^{1/q},R) \cong H^{0}(X, \O_{X}((1-q)K_{X}))^{1/q}$ derived from the adjunction formula, we may regard $\phi$ and $\psi$ as rational sections of the sheaf $\omega_{\widetilde{X}}^{(1-q)}$ and consider the corresponding divisors on $\widetilde{X}$
$$D_{\phi}=(\phi)_0-(\phi)_{\infty}, \quad D_{\psi}=(\psi)_0-(\psi)_{\infty},$$
where $(\phi)_0$ and $(\phi)_{\infty}$ (resp. $(\psi)_0$ and $(\psi)_{\infty}$) are the divisors of zeros and poles of $\phi$ (resp. $\psi$) as a rational section of $\omega_{\widetilde{X}}^{(1-q)}$.
Clearly $D_{\phi}=D_{\psi}+\Div_{\widetilde{X}}(c)+\Div_{\widetilde{X}}(d)$.
By definition, $D_{\phi}$ and $D_{\psi}$ are linearly equivalent to $(1-q)K_{\widetilde{X}}$, and $(\phi)_{\infty}$ and $(\psi)_{\infty}$ are $f$-exceptional divisors.
We denote $X'=\widetilde{X} \setminus \Supp (\phi)_{\infty}$. Then 
$$\phi \in \Hom_{\O_{X'}}(\O_{X'}^{1/q}, \O_{X'}) \cong H^0(X', \O_{X'}((1-q)K_{X'}))^{1/q}.$$
\begin{claim}
The coefficient of $D_{\phi}$ in each irreducible component is less than or equal to $q-1$.
\end{claim}
\begin{proof}[Proof of Claim]
Assume to the contrary that there exists an irreducible component $D_{\phi,0}$ of $D_{\phi}$ whose coefficient is greater than $q-1$.
Then $\phi$ lies in 
$$\Hom_{\O_{X'}}(\O_{X'}(qD_{\phi,0})^{1/q}, \O_{X'}) \cong H^0(X', \O_{X'}((1-q)K_{X'}-qD_{\phi,0}))^{1/q}$$
and gives a splitting of the map $\O_{X'} \hookrightarrow \O_{X'}(qD_{\phi,0})^{1/q}$. This map, however, factors through $\O_{X'}(D_{\phi,0})$ and $\O_{X'} \hookrightarrow \O_{X'}(D_{\phi,0})$ never splits as an $\O_{X'}$-module homomorphism. This is a contradiction.
\end{proof}

Let $B=\frac{1}{q-1}D_{\psi}$ and 
then $B$ is $\Q$-linearly equivalent to $-K_{\widetilde{X}}$, so that $f_*B$ is $\Q$-linearly equivalent to $-K_X$.
Hence $f_*B$ is an effective $\Q$-Cartier divisor and $(B-f^*f_*B)+ \sum_{j=1}^n a_jE_j+\sum_{i=1}^k t_iG_i$ is an $f$-exceptional divisor which is $\Q$-linearly trivial relative to $f$. Hence 
$$(B-f^*f_*B)+\sum_{j=1}^n a_jE_j+\sum_{i=1}^k t_iG_i=0.$$

Since $D_{\phi}=D_{\psi}+\Div_{\widetilde{X}}(c)+\Div_{\widetilde{X}}(d)$ and $\Div_{\widetilde{X}}(d) \ge \sum_{i=1}^k\lfloor t_i(q-1) \rfloor G_i$, by the above claim, the coefficient of $D_{\psi}+\sum_{i=1}^k\lfloor t_i(q-1) \rfloor G_i$ in $E_j$ is less than or equal to $q-1$ for $1 \le j \le m$ and the coefficient of $D_{\psi}+\sum_{i=1}^k t_i(q-1) G_i$ in $E_j$ is less than $q-1$ for $m+1 \le j \le n$.
Since 
\begin{align*}
B-f^*f_*B +\sum_{i=1}^k t_iG_i & \le \frac{1}{q-1}(D_{\psi}+\sum_{i=1}^k t_i(q-1) G_i)\\
& < \frac{1}{q-1}(D_{\psi}+\sum_{i=1}^k\lfloor t_i(q-1) \rfloor G_i)+\frac{\sum_{i=1}^kG_i}{q-1},
\end{align*}
the coefficient of $B-f^*f_*B+\sum_{i=1}^k t_iG_i$ in $E_j$ is less than $1+\frac{M}{q-1}$ for $1 \le j \le m$ and is less than $1$ for $m+1 \le j \le n$, where $M$ is a constant and independent of $q$.
By taking $q=p^e$ sufficiently large, we see that the coefficient of $B-f^*f_*B+\sum_{i=1}^k t_iG_i$ in $E_j$ is less than or equal to $1$ for $1 \le j \le m$. Thus we have $a_j \ge -1$ for every $j=1,\dots, m$ and $a_j > -1$ for every $j=m+1, \dots, n$, which implies the pair $(X,Y)$ is plt.
\end{proof}

\begin{lem}\label{Fedder}
{\rm (Fedder type criteria)}
Let $(R,\m)$ be an F-finite regular local ring of characteristic $p>0$ and $I \subsetneq R$ a reduced ideal. Let $\a_1, \dots, \a_k$ be ideals of $R$ and $t_1, \dots, t_k > 0$ real numbers. Write $S=R/I$.
\begin{enumerate}
\item
The pair $(S, (\a_1S)^{t_1}\cdots(\a_kS)^{t_k})$ is F-pure if and only if for all large $q=p^e$, $\a_1^{\lfloor t_1(q-1) \rfloor}\cdots \a_k^{\lfloor t_k(q-1) \rfloor}(I^{[q]}:I) \not\subset \m^{[q]}$.
\item
The pair $(S, (\a_1S)^{t_1}\cdots(\a_kS)^{t_k})$ is purely F-regular if and only if for every element $c \in R \setminus I$ which is not in any minimal prime ideal of $\a_1 \cdots \a_k+I$, there exists $q'$ such that for every $q=p^e \ge q'$, $c\a_1^{\lfloor t_1(q-1) \rfloor}\cdots \a_k^{\lfloor t_k(q-1) \rfloor}(I^{[q]}:I) \not\subset \m^{[q]}$.
\item
The pair $(S, (\a_1S)^{t_1}\cdots(\a_kS)^{t_k})$ is strongly F-regular if and only if for every element $c \in R \setminus I$, there exists $q=p^e$ such that $c\a_1^{\lfloor t_1(q-1) \rfloor}\cdots \a_k^{\lfloor t_k(q-1) \rfloor}(I^{[q]}:I) \not\subset \m^{[q]}$.
By Proposition \ref{basic} (3), this is equivalent to saying that for every element $c \in R \setminus I$, there exists $q'$ such that for every $q=p^e \ge q'$, $c\a_1^{\lceil t_1q \rceil}\cdots \a_k^{\lceil t_kq \rceil}(I^{[q]}:I) \not\subset \m^{[q]}$.
\end{enumerate}
\end{lem}
\begin{proof}
The proof is similar to those in \cite{Fe}, \cite{Gl} and \cite{HW}. 
We may assume without loss of generality that $(R,\m)$ is a $d$-dimensional complete regular local ring. Let $E_R$ and $E_S$ be injective hulls of the residue fields of $R$ and $S$ respectively. Since $R$ is regular, in particular Gorenstein, $E_R \cong H^d_{\m}(R)$ and $E_R \otimes_R R^{1/q} \cong H^d_{\m}(R^{1/q})$.
We can identify $E_R$ with $E_R \otimes_R R^{1/q}$ via the identification of $R$ with $R^{1/q}$ and view $E_S$ as a submodule of $E_R$ via the isomorphism $E_S \cong (0:I)_{E_R} \subset E_R$. $E_R$ and $E_S$ have the one-dimensional socle in common, and let $z \in E_R$ a generator of the socle.
Since $R \hookrightarrow R^{1/q}$ is flat by Lemma \ref{flat}, via the identification $R \cong R^{1/q}$, we have $E_S \otimes_R R^{1/q} \cong (0:I^{[q]})_{E_R}$ in $E_R \otimes_R R^{1/q} \cong E_R$. Therefore
$$
\begin{small}
E_S \otimes_S S^{1/q} \cong E_S \otimes_R R^{1/q} \otimes_{R^{1/q}} S^{1/q} \cong (0:I^{[q]})_{E_R} \otimes_R S \cong \frac{(0:I^{[q]})_{E_R}}{I(0:I^{[q]})_{E_R}}.
\end{small}
$$
Here notice that by the Matlis duality (cf. \cite[Lemma 3.3]{Ha3}), 
$I(0:I^{[q]})_{E_R}=\Ann_{E_R} \Ann_R I(0:I^{[q]})_{E_R}=(0:(I^{[q]}:I))_{E_R}$, because $R$ is complete.

Let $d \in R$ be any nonzero element and $F^e_R:E_R \to E_R \otimes_R R^{1/q} \cong E_R$ (resp. $F^e_S:E_S \to E_S \otimes_S S^{1/q}$) the $e$-times iterated Frobenius map induced on $E_R$ (resp. $E_S$). 
By the above argument, $dF^e_S:E_S \to E_S \otimes_S S^{1/q}$ is injective if and only if $dF^e_S(z) \ne 0$ if and only if $dF^e_R(z) \notin (0:(I^{[q]}:I))_{E_R}$. Since $F^e_R(z) \in E_R$ generates $(0:\m^{[q]})_{E_R}$, this is equivalent to saying that $d(0:\m^{[q]})_{E_R} \not\subset (0:(I^{[q]}:I))_{E_R}$, namely $d(I^{[q]}:I) \not\subset \m^{[q]}$. Thus, by Lemma \ref{injective}, the pair $(S, (\a_1S)^{t_1}\cdots(\a_kS)^{t_k})$ is F-pure if and only if for every large $q=p^e$, $\a_1^{\lfloor t_1(q-1) \rfloor}\cdots \a_k^{\lfloor t_k(q-1) \rfloor}(I^{[q]}:I) \not\subset \m^{[q]}$. $(2)$ and $(3)$ are derived from a similar argument.
\end{proof}

\begin{rem}\label{Fedder2}
Thanks to Corollary \ref{F-regular} and Lemma \ref{test}, Lemma \ref{Fedder} $(3)$ has a few variants. Assume that one of the following conditions is satisfied. 
\renewcommand{\labelenumi}{(\alph{enumi})}
\begin{enumerate}
\item
$S$ is $\Q$-Gorenstein and the image of an element $c \in R \setminus I$ in $S$ is an $(\a_1S)^{t_1}\cdots(\a_kS)^{t_k}$-test element. 
\item
$c \in R \setminus I$ is an element such that the localization $S_c=R_c/IR_c$ with respect to $c$ is strongly F-regular. 
\end{enumerate}
Then the pair $(S, (\a_1S)^{t_1} \cdots (\a_kS)^{t_k})$ is strongly F-regular if and only if there exists $q=p^e$ such that $c\a_1^{\lceil t_1q \rceil}\cdots \a_k^{\lceil t_kq \rceil}(I^{[q]}:I) \not\subset \m^{[q]}$. 
\end{rem}

The following theorem is a generalization of ``F-inversion of Adjunction'' \cite[Theorem 4.9]{HW}.
\begin{thm}\label{F-inversion}
Let $(R, \m)$ be an F-finite regular local ring of characteristic $p>0$ and $I \subsetneq R$ an unmixed reduced ideal of height $h>0$. Let $\a_1, \dots, \a_k$ be ideals of $R$ and $t_1, \dots, t_k > 0$ real numbers. We denote $S=R/I$.
\begin{enumerate}
\item
If the pair $(S, (\a_1S)^{t_1}\cdots (\a_kS)^{t_k})$ is strongly F-regular $($resp. F-pure$)$, then the pair $(R, I\a_1^{t_1}\cdots \a_k^{t_k})$ is purely F-regular $($resp. F-pure$)$. 
\item
If the pair $(R, I^h\a_1^{t_1}\cdots \a_k^{t_k})$ is F-pure, then $(S, (\a_1S)^{t_1}\cdots (\a_kS)^{t_k})$ is an F-pure pair. 
In addition, suppose all $\a_i$ contain $I$. 
If the pair $(R, I^h\a_1^{t_1}\cdots \a_k^{t_k})$ is purely F-regular, then $(S, (\a_1S)^{t_1}\cdots (\a_kS)^{t_k})$ is a strongly F-regular pair.
\end{enumerate}
\end{thm}
\begin{proof}
$(1)$ 
First we consider the case where $(S, (\a_1S)^{t_1}\cdots (\a_kS)^{t_k})$ is an F-pure pair. 
By Lemma \ref{Fedder}, the pair $(S, (\a_1S)^{t_1}\cdots (\a_kS)^{t_k})$ (resp. $(R, I\a_1^{t_1}\cdots \a_k^{t_k})$) is F-pure if and only if $\a_1^{\lfloor t_1(q-1) \rfloor}\cdots \a_k^{\lfloor t_k(q-1) \rfloor}(I^{[q]}:I) \not\subset \m^{[q]}$ (resp. $\a_1^{\lfloor t_1(q-1) \rfloor}\cdots \a_k^{\lfloor t_k(q-1) \rfloor}I^{q-1} \not\subset \m^{[q]}$) for all large $q=p^e$. 
Therefore it suffices to show that $(I^{[q]}:I) \subset I^{q-1}$ for every $q=p^e$.
Fix any element $u \in (I^{[q]}:I)$ and any power $Q=p^l \ge q=p^e$.
Write $Q=a(q-1)+r$ with $0 \le r < q-1$.
Then $a=p^{l-e}+p^{l-2e}+\dots+p^{l-en}$ and $r=p^{l-en}$ where $n$ is the maximal integer which is not greater than $l/e$.
Let $I=\cap_{i=1}^m P_i$ be the irredundant prime decomposition of $I$, that is, $P_i$ are minimal prime ideals of $I$. 
Note that, by \cite[Lemma 4.1]{Fe}, $(I^{[q]}:I)=\cap_{i=1}^m (P_i^{[q]}:P_i)$.

\begin{claim}
For all $i=1,\dots,m$,
$$I^{h(p^{e-1}-1)+1}u^a \subset (P_iR_{P_i})^{[Q]}.$$
\end{claim}

\begin{proof}[Proof of Claim]
Let $x_{i,1}, \dots, x_{i,h}$ be a regular system of parameters of $R_{P_i}$, and denote $x_i=x_{i,1} \cdots x_{i,h} \in R_{P_i}$.
Then 
\begin{align*}
(P_i^{[q]}:P_i)R_{P_i}&=((P_iR_{P_i})^{[q]}:P_iR_{P_i})\\
&={x_i}^{q-1}R_{P_i}+({x_{i,1}}^q, \dots, {x_{i,h}}^q)R_{P_i}.
\end{align*}
Therefore we can write in $R_{P_i}$
$$u=a_{i,0}{x_i}^{q-1}+\sum_{j=1}^h a_{i,j} {x_{i,j}}^q,$$ where $a_{i,j} \in R_{P_i}$ for every $j=0, 1, \dots, h$. Then
\begin{align*}
u^a&=\prod_{k=1}^n u^{p^{l-ke}}\\
&=\prod_{k=1}^n \left({a_{i,0}}^{p^{l-ke}}{x_i}^{p^{l-(k-1)e}-p^{l-ke}}+\sum_{j=1}^h {a_{i,j}}^{p^{l-ke}} {x_{i,j}}^{p^{l-(k-1)e}}\right)\\
& \in {x_i}^{Q-r}R_{P_i}+(P_iR_{P_i})^{[Q]}.
\end{align*}
Since 
$$I^{h(p^{e-1}-1)+1}R_{P_i} \subset (P_iR_{P_i})^{h(r-1)+1} \subset ({x_{i,1}}^r, \ldots, {x_{i,h}}^r)R_{P_i},$$ 
we obtain the inclusion
$I^{h(p^{e-1}-1)+1}u^a \subset (P_iR_{P_i})^{[Q]}$.
\end{proof}

Since $P_i^{[Q]}=(P_iR_{P_i})^{[Q]} \cap R$ by Lemma \ref{flat}, it follows from the above claim that $I^{h(p^{e-1}-1)+1}u^a \subset P_i^{[Q]}$ for every $i=1, \dots, m$. By Lemma \ref{flat} again, $I^{h(p^{e-1}-1)+1}u^a \subset \cap_{i=1}^m P_i^{[Q]}=I^{[Q]}$. Taking the $(q-1)$-th powers of both sides and abbreviating $b=\{h(p^{e-1}-1)+1\}(q-1)$, we have $I^{b}u^{a(q-1)} \subset (I^{q-1})^{[Q]}$, in particular $I^{b}u^{Q} \subset (I^{q-1})^{[Q]}$ for every $Q \ge q$. Notice that $b$ does not depend on $Q$. Thus $u \in (I^{q-1})^*=I^{q-1}$, because every ideal is tightly closed in case the ring is regular (see \cite[Theorem 4.6]{HH1}). 

Next we suppose that the pair $(S, (\a_1S)^{t_1}\cdots (\a_kS)^{t_k})$ is strongly F-regular. We may assume that $I$ does not contain $\a_i$ for any $i=1, \dots, k$. 
Thanks to Lemma \ref{Fedder}, the pair $(S, (\a_1S)^{t_1}\cdots (\a_kS)^{t_k})$ (resp. $(R, I\a_1^{t_1}\cdots \a_k^{t_k})$) is strongly F-regular (resp. purely F-regular) if and only if for every element $c \in R \setminus I$ (resp. for every element $c \in R^{\circ}$ which is not in any minimal prime ideal of $I\a_1 \cdots \a_k$), there exists $q'$ such that for every $q=p^e \ge q'$, $c\a_1^{\lceil t_1q \rceil}\cdots \a_k^{\lceil t_kq \rceil}(I^{[q]}:I) \not\subset \m^{[q]}$ (resp. $c\a_1^{\lfloor t_1(q-1) \rfloor}\cdots \a_k^{\lfloor t_k(q-1) \rfloor}I^{q-1} \not\subset \m^{[q]}$). 
Since $S$ is strongly F-regular by Proposition \ref{basic} $(1)$, $I$ is a minimal prime ideal of $I\a_1 \cdots \a_k$. 
Thus it is enough to show that $(I^{[q]}:I) \subset I^{q-1}$ for every $q=p^e$, but it has been already proved in the F-pure case. 

$(2)$ 
First we consider the F-pure case. 
Thanks to Lemma \ref{Fedder}, it is enough to prove that $I^{h(q-1)} \subset (I^{[q]}:I)$ for all $q=p^e$. Let $I=\cap_{i=1}^m P_i$ be the irredundant prime decomposition of $I$. By Lemma \ref{flat}, an element $x \in R^{\circ}$ is a nonzero divisor on $R/P_i$ if and only if $x$ is a nonzero divisor on $R/P_i^{[q]}$. Hence, if $x$ is a nonzero divisor on $R/P_i$, then $x$ is a nonzero divisor on $R/(P_i^{[q]}:P_i)$. Thus $(P_i^{[q]}:P_i)$ is $P_i$-primary. Since $R_{P_i}$ is a regular local ring of dimension $h$, 
\begin{align*}
(P_i^{[q]}:P_i)&=((P_iR_{P_i})^{[q]}:P_iR_{P_i}) \cap R\\
&=(P_i^{h(q-1)}R_{P_i}+P_i^{[q]}R_{P_i}) \cap R \supset P_i^{h(q-1)}.
\end{align*}
Applying \cite[Lemma 4.1]{Fe}, we have
$$I^{h(q-1)} \subset \displaystyle\bigcap_{i=1}^m P_i^{h(q-1)} \subset \displaystyle\bigcap_{i=1}^m (P_i^{[q]}:P_i) = (I^{[q]}:I).$$

Next we suppose that the pair $(R, I^h\a_1^{t_1} \cdots \a_k^{t_k})$ is purely F-regular.  
Since every $\a_i$ contains $I$ by assumption, then all the minimal prime ideals of $I\a_1 \cdots \a_k$ are minimal prime ideals of $I$. 
Therefore, by Lemma \ref{Fedder}, it is sufficient to show that $I^{h(q-1)} \subset (I^{[q]}:I)$ for all $q=p^e$, but it has been already proved in the F-pure case. 
\end{proof}

\begin{rem}
$(1)$
In \cite[Theorem 4.9]{HY}, the assumption that $R$ is regular is unnecessary.
However, when the height of $I$ is greater than one, we cannot drop the regularity of $R$ from the statement of Theorem \ref{F-inversion}. See Example \ref{example} (i).

$(2)$
If the ideal $I$ is generated by a regular sequence, then it follows from the repeated applications of ``F-Inversion of Adjunction'' \cite[Theorem 4.9]{HY} that if the pair $(S, (\a_1S)^{t_1}\cdots (\a_kS)^{t_k})$ is F-pure (resp. strongly F-regular), then $(R, I^h\a_1^{t_1}\cdots \a_k^{t_k})$ is an F-pure pair (resp. a purely F-regular pair).
However, when $I$ is not generated by a regular sequence, the pair $(R, I^h\a_1^{t_1}\cdots \a_k^{t_k})$ is not necessarily F-pure (resp. purely F-regular), even if $(S, (\a_1S)^{t_1}\cdots (\a_kS)^{t_k})$ is F-pure (resp. strongly F-regular). See Example \ref{example} (ii), (iii). 
\end{rem}

\begin{cor}\label{restriction}
Let $(R,\m)$ be an F-finite regular local ring of characteristic $p>0$ and $I \subset R$ a nonzero prime ideal such that $R/I$ is a $\Q$-Gorenstein normal local ring. Let $\a_1, \dots, \a_k$ be ideals of $R$ and $t_1, \dots, t_k > 0$ real numbers. Then, setting $S=R/I$, we have
$$\tau(S, (\a_1S)^{t_1}\cdots (\a_kS)^{t_k}) \subset \tau(R, I^t\a_1^{t_1}\cdots \a_k^{t_k})S$$
for every $0 < t < 1$.
\end{cor}
\begin{proof}
By \cite[Proposition 3.2]{HT}, we may assume without loss of generality that $R$ is complete. 
Let $E_R=E_R(R/\m)$ (resp. $E_S=E_S(S/\m S)$) be the injective hull of the residue field of $R$ (resp. $S$). For $q=p^e$, let $F^e_R:E_R \to E_R \otimes R^{1/q} \cong E_R$ (resp. $F^e_S:E_S \to E_S \otimes S^{1/q}$) be the $e$-times iterated Frobenius map induced on $E_R$ (resp. $E_S$). 
We can view $E_S$ as a submodule of $E_R$ via the isomorphism $E_S \cong (0:I)_{E_R} \subset E_R$.
\begin{claim}
$$0_{E_R}^{*I^t\a_1^{t_1}\cdots \a_k^{t_k}} \cap E_S \subset 0_{E_S}^{*(\a_1S)^{t_1}\cdots (\a_kS)^{t_k}}.$$
\end{claim}
\begin{proof}[Proof of Claim]
Let $z \in 0_{E_R}^{*I^t\a_1^{t_1}\cdots \a_k^{t_k}} \cap E_S$. 
Since the unit element $1$ is an $I^t\a_1^{t_1}\cdots \a_k^{t_k}$-test element by Theorem \ref{exist}, for all $q=p^e$,
$$I^{\lceil tq \rceil}\a_1^{\lceil t_1q \rceil}\cdots \a_k^{\lceil t_kq \rceil}F^e_R(z) =0 \in E_R \otimes_R R^{1/q} \cong E_R.$$
We choose a sufficiently large $q=p^e$ so that $\lceil tq \rceil \le q-1$.
By the proof of Theorem \ref{F-inversion}, we have $(I^{[q]}:I) \subset I^{q-1} \subset I^{\lceil tq \rceil}$. Hence
$$\a_1^{\lceil t_1q \rceil}\cdots \a_k^{\lceil t_kq \rceil}F^e_R(z) \in (0:(I^{[q]}:I))_{E_R}.$$
By the same argument as that of Lemma \ref{Fedder}, for all large $q=p^e$,
$$ \a_1^{\lceil t_1q \rceil}\cdots \a_k^{\lceil t_kq \rceil}F^e_S(z)=0 \in E_S \otimes_S S^{1/q},$$
whence $z$ is in $0_{E_S}^{*(\a_1S)^{t_1}\cdots (\a_kS)^{t_k}}$.
\end{proof}
Since $R$ is complete, we have $0_{E_R}^{*I^t\a_1^{t_1}\cdots \a_k^{t_k}}=(0:\tau(R, I^t\a_1^{t_1}\cdots \a_k^{t_k}))_{E_R}$. Hence
\begin{align*}
0_{E_R}^{*I^t\a_1^{t_1}\cdots \a_k^{t_k}} \cap E_S&=(0:\tau(R, I^t\a_1^{t_1}\cdots \a_k^{t_k})+I)_{E_R}\\
&=\Bigl(0:\frac{\tau(R, I^t\a_1^{t_1}\cdots \a_k^{t_k})+I}{I}\Bigr)_{E_S}\\
&=(0:\tau(R, I^t\a_1^{t_1}\cdots \a_k^{t_k})S)_{E_S}.
\end{align*}
Since $S$ is $\Q$-Gorenstein, by Definition-Theorem \ref{taudef} and the above claim, we have 
\begin{align*}
\tau(S, (\a_1S)^{t_1}\cdots (\a_kS)^{t_k})&=\Ann_S(0_{E_S}^{*(\a_1S)^{t_1}\cdots (\a_kS)^{t_k}})\\
 &\subset \Ann_S({0_{E_R}^{*I^t\a_1^{t_1}\cdots \a_k^{t_k}} \cap E_S})\\
&=\tau(R, I^t\a_1^{t_1}\cdots \a_k^{t_k})S.
\end{align*} 
\end{proof}

\begin{expl}\label{example}
For a pair $(R, \a)$ of an F-finite strongly F-regular ring $R$ of characteristic $p>0$ and an ideal $\a \subset R$, we denote by $c(R,\a)$ the {\it F-pure threshold} of $\a$, that is,
\begin{align*}
c(R,\a)&=\mathrm{Sup}\{t \in \R_{\ge 0} \mid (R, \a^t) \text{ is F-pure} \}\\
&=\mathrm{Sup}\{t \in \R_{\ge 0} \mid (R, \a^t) \text{ is strongly F-regular} \}.
\end{align*}
See \cite{TW} for properties and computations of F-pure thresholds.
\renewcommand{\labelenumi}{(\roman{enumi})}
\begin{enumerate}
\item Let $R=k[[X, Y, Z]]/(X^2+Y^3+Z^5)$ be the rational double point of type $E_8$ over a field $k$ of characteristic $p>5$ and $\m$ the maximal ideal of $R$. 
By Lemma \ref{Fedder}, the pair $(R, \m^t)$ is F-pure if and only if for all large $q=p^e$, $(X^2+Y^3+Z^5)^{q-1}(X,Y,Z)^{\lfloor t(q-1) \rfloor} \not\subset (X^q, Y^q, Z^q)$ in $k[[X,Y,Z]]$.  
Hence the F-pure threshold $c(R,\m)$ is bounded by the supremum of all real numbers $t>0$ such that for all large $q=p^e$, 
$$ X^{2 \cdot \frac{q-1}{2}}Y^{3\lfloor \frac{q-1}{3} \rfloor}Z^{5(q-1-\frac{q-1}{2}-\lfloor \frac{q-1}{3} \rfloor)}(X,Y,Z)^{\lfloor t(q-1) \rfloor} \not\subset (X^q,Y^q,Z^q). $$
Thus we have $c(R,\m) \le 1/6$, particularly $R$ is not F-pure with respect to $\m$, while $R/\m$ is a field and in particular strongly F-regular. 
\item Let $R=k[X_1, \dots, X_d]$ be a polynomial ring of dimension $d \ge 3$ over a field $k$ of characteristic $p>0$ and 
$$I=\sum_{i=1}^d (X_1 \cdots X_{i-1} X_{i+1} \cdots X_d) R \subset R.$$ 
Then $I$ is an ideal of height $d-1$ and $R/I$ is F-pure. 
On the other hand, since the Newton polytope $P(I)$ associated to $I$ is $\{u_1+ \dots +u_d \ge d-1\}$ ($u_1, \dots, u_d$ are the natural coordinates on $\R^d$), by Lemma \ref{F-regular} and \cite[Theorem 4.6]{HY}, $c(R,I)=\frac{d}{d-1}$ and in particular the pair $(R,I^{d-1})$ is not F-pure. 
\item 
Let $R=k[X^aY^bZ^c\mid a+b+c \equiv 0 \mod 3] \subset k[X,Y,Z]$ be the cyclic quotient singularity of type 1/3(1,1,1) over a field $k$ of characteristic $p>0$. 
Let $S=k[T_1, \dots, T_{10}]$ be a ten-dimensional polynomial ring over $k$ and $\m_S=(T_1, \dots, T_{10}) \subset S$ a maximal ideal. 
$\Spec R$ is naturally embedded in $\mathbb{A}_k^{10}=\Spec S$ and denote by $I \subset S$ the defining ideal of $\Spec R$ in $\mathbb{A}_k^{10}$. 
Then $I$ is an ideal of height $7$ generated by quadratics. 
Since $I^{5q} \subset \m_S^{10q} \subset \m_S^{[q]}$ for all $q=p^e$, by Lemma \ref{Fedder}, the pair $(S, I^5)$ is not strongly F-regular, that is, $c(S,I) \le 5$. Thus the pair $(S, I^7)$ is not F-pure, while $S/I$ is toric and in particular strongly F-regular (or klt). 
\end{enumerate}
\end{expl}

\section{Inversion of Adjunction of arbitrary codimension}
As an application of Theorem \ref{F-inversion}, we prove a sort of Inversion of Adjunction  of arbitrary codimension for klt pairs.
\begin{thm}\label{klt}
Let $X$ be a non-singular variety over a field of characteristic zero, $\Delta$ an effective $\R$-divisor on $X$ and $Y=\sum_{i=1}^k t_iY_i$ a formal combination where $t_i > 0$ are real numbers and $Y_i \subsetneq X$ are closed subschemes. Let $Z \subsetneq X$ be a normal closed subvariety such that $Z \not\subset \Delta$ and $Z \not\subset \cup_{i=1}^{k}Y_i$. 
\begin{enumerate}
\item If the pair $(Z, \Delta|_Z)$ is admissible $($see  Definition \ref{admissible}$)$, then the pair $(X,\Delta+Z)$ is plt near $Z$. 
\item Suppose $Z$ is $\Q$-Gorenstein. 
If the pair $(Z, Y|_Z)$ is klt, then the pair $(X, Y+Z)$ is plt near $Z$.
\end{enumerate}
\end{thm}
\begin{proof}
The questions are local, so we may assume that $X=\Spec R$ and $Z=\Spec S$, where $R$ is a regular local ring essentially of finite type over a field of characteristic zero and $S=R/I$ with $I$ a nonzero ideal of $R$.

$(1)$
By Lemma \ref{admissible}, the pair $(Z, \Delta|_Z)$ is of strongly F-regular type.
Since $\Delta$ is an effective $\R$-Cartier divisor, we can write $\Delta=\sum_{j=1}^l s_j \cdot \Div(f_j)$ with $f_j \in R^{\circ}$ and $s_j \in \R_{> 0}$ for all $j=1, \dots, l$. By Remark \ref{divisor}, the strong F-regularity of the pair $(Z, \Delta|_Z)$ is equivalent to the strong F-regularity of the pair $(S, (f_1S)^{s_1} \cdots (f_lS)^{s_l})$, hence it follows from Theorem \ref{F-inversion} that the pair $(R, f_1^{s_1} \cdots f_l^{s_l}I)$ is of purely F-regular type.
By Proposition \ref{plt}, the pair $(X, \Delta+Z)$ is plt near $Z$.

$(2)$ 
Let $\a_i \subset R$ be the defining ideal of the closed subscheme $Y_i \subsetneq X$ for all $i=1, \dots, n$.
By virtue of Theorem \ref{multiplier} and Corollary \ref{F-regular}, the pair $(S, (\a_1S)^{t_1} \cdots (\a_kS)^{t_k})$ is of strongly F-regular type. It follows from Theorem \ref{F-inversion} that $(R,I\a_1^{t_1} \cdots \a_k^{t_k})$ is a pair of purely F-regular type, hence the pair $(X,Y+Z)$ is plt near $Z$ by Proposition \ref{plt}. 
\end{proof}

Using Corollary \ref{restriction}, we show the lc case of Theorem \ref{klt} $(2)$. 
The case where $Z$ is a divisor was proved by Ein, Musta{\c{t}}{\v{a}} and Yasuda \cite[Corollary 1.8]{EMY} and we give another proof of their result by characteristic $p$ methods.
\begin{thm}\label{inversion}
Let $X$ be a non-singular variety over a field of characteristic zero and $Y=\sum_{i=1}^k t_iY_i$ a formal combination where $t_i > 0$ are real numbers and $Y_i \subsetneq X$ are closed subschemes. Let $Z \subsetneq X$ be a normal $\Q$-Gorenstein closed subvariety such that $Z \not\subset \cup_{i=1}^k Y_i$. 
If the pair $(Z, Y|_Z)$ is lc, then the pair $(X,Y+Z)$ is also lc near $Z$.
\end{thm}
\begin{proof}
The question is local, so we may assume that $X=\Spec R$ and $Z=\Spec S$, where $R$ is a regular local ring essentially of finite type over a field of characteristic zero with the maximal ideal $\m$ and $S=R/I$ with $I \subset R$ an ideal of positive height.
Let $J \subseteq R$ be the reduced ideal containing $I$ which defines the non-klt locus of $Z$ and we denote by $\a_i \subset R$ the nonzero ideal associated to the closed subscheme $Y_i$ for all $i=1,2, \dots, n$.
Since the pair $(Z, Y|_Z)$ is lc, thanks to Lemma \ref{lc}, we have 
$$J\O_Z \subset \J(Z,(\a_1\O_Z)^{st_1}\cdots (\a_k\O_Z)^{st_k}(J\O_Z)^s)$$ 
for every $0<s<1$.  Now fix $0<s<1$ and we consider the reduction from characteristic zero to characteristic $p \gg 0$. Then it follows from Theorem \ref{multiplier} that 
$$JS \subset \tau(S, (\a_1S)^{st_1}\cdots (\a_kS)^{st_k}(JS)^s).$$
Then the following claim is essential. 

\begin{claim}
$$\tau(R, \a_1^{st_1}\cdots \a_k^{st_k}J^s)=R.$$
\end{claim}

\begin{proof}[Proof of Claim]
We may assume without loss of generality that $R$ is complete. 
Let $E_R=E_R(R/\m)$ (resp. $E_S=E_S(S/\m S)$) be the injective hull of the residue field of $R$ (resp. $S$). For $q=p^e$, let $F^e_R:E_R \to E_R \otimes R^{1/q}$ (resp. $F^e_S:E_S \to E_S \otimes S^{1/q}$) be the $e$-times iterated Frobenius map induced on $E_R$ (resp. $E_S$).
Since $JS \subset \tau(S, (\a_1S)^{st_1}\cdots (\a_kS)^{st_k}(JS)^s)$, $J 0_{E_S}^{*(\a_1S)^{st_1}\cdots (\a_kS)^{st_k}(JS)^s}=0$ by definition.
We identify $E_S$ with $(0:I)_{E_R} \subset E_R$ and fix an element $y \in (0:J\m)_{E_S} \setminus (0:J)_{E_S}$.
Then $y$ does not belong to $0_{E_S}^{*(\a_1S)^{st_1}\cdots (\a_kS)^{st_k}(JS)^s}$, that is, for every $c \in R \setminus I$,
there exists $q=p^e$ such that
$$c\a_1^{\lceil st_1q \rceil} \cdots \a_k^{\lceil st_kq \rceil} J^{\lceil sq \rceil}F^e_S(y) \ne 0 \in E_S \otimes_S S^{1/q}.$$
By a similar argument to that of Lemma \ref{Fedder}, 
$$c\a_1^{\lceil st_1q \rceil} \cdots \a_k^{\lceil st_kq \rceil}J^{\lceil sq \rceil}F^e_R(y) \not\subset I(0:I^{[q]})_{E_R}=(0:(I^{[q]}:I))_{E_R}.$$
Since $F^e_R(y) \in (0:(J\m)^{[q]})_{E_R}$, we have
\begin{align}
c\a_1^{\lceil st_1q \rceil} \cdots \a_k^{\lceil st_kq \rceil}J^{\lceil sq \rceil}(I^{[q]}:I) \not\subset (J\m)^{[q]}.
\end{align}

On the other hand, for every minimal prime ideal $P$ of $J$, by the definition of the ideal $J$, $\Spec S_P$ is not klt, in particular not strongly F-regular (cf. \cite[Theorem 3.3]{HW}). Therefore, by Lemma \ref{Fedder}, there exists some $c_P \in R \setminus I$ such that for every power $Q$ of $p$,
$$c_P(I^{[Q]}:I)R_P=c_P((IR_P)^{[Q]}:IR_P) \subset (PR_P)^{[Q]}.$$
Since $P^{[Q]}=(PR_P)^{[Q]} \cap R$ by Lemma \ref{flat},
$c_P(I^{[Q]}:I) \subset P^{[Q]}$ for every minimal prime ideal $P$ of $J$.
Since $J$ is reduced, by the flatness of the Frobenius map again, there exists some $c \in R \setminus I$ such that 
\begin{align}
c(I^{[Q]}:I) \subset \bigcap_{P \supset J} P^{[Q]}=J^{[Q]}
\end{align}
for every power $Q$ of $p$, where $P$ runs through all minimal prime ideals of $J$.
Thus comparing $(4.1)$ with $(4.2)$, we have
$$\a_1^{\lceil st_1q \rceil} \cdots \a_k^{\lceil st_kq \rceil}J^{\lceil sq \rceil} \not\subset \m^{[q]}$$
for some $q=p^e$.
Since the unit element $1$ is an $\a_1^{st_1} \cdots \a_k^{st_k}J^{s}$-test element by Theorem \ref{exist},
it follows from Remark \ref{Fedder2} that the pair $(R, \a_1^{st_1}\cdots \a_k^{st_k}J^s)$ is strongly F-regular, namely $\tau(R, \a_1^{st_1}\cdots \a_k^{st_k}J^s)=R$. 
\end{proof}

It follows from Corollary \ref{restriction} that
$$JS \subset \tau(S, (\a_1S)^{st_1}\cdots (\a_kS)^{st_k}(JS)^s) \subset \tau(R, \a_1^{st_1}\cdots \a_k^{st_k}I^sJ^s)S. $$
Hence $J \subset \tau(R, \a_1^{st_1}\cdots \a_k^{st_k}I^sJ^s)$, because $J$ contains $I$ and $I \subset \tau(R, \a_1^{st_1}\cdots \a_k^{st_k}I^sJ^s)$ by the above claim.
Now we return to the situation in characteristic zero and we have 
$$J \subset \J(X, \a_1^{st_1}\cdots \a_k^{st_k}I^sJ^s)$$
for all $0<s<1$. By Lemma \ref{lc}, this implies the pair $(X, Y+Z)$ is lc near $Z$.
\end{proof}

\begin{rem}
Ein and Musta{\c{t}}{\v{a}} \cite[Corollary 3.2]{EM} proved that in case $Z$ is a locally complete intersection variety of codimension $r>0$, the pair $(Z, Y|_Z)$ is lc if and only if the pair $(X,Y+rZ)$ is also lc near $Z$. However, when $Z$ is not locally complete intersection, the pair $(X,Y+rZ)$ is not necessarily lc near $Z$, even if the pair $(Z, Y|_Z)$ is lc. (cf. Example \ref{example}).
\end{rem}

Finally we remark on a sort of Adjunction of arbitrary codimension. 
We expect that pairs of dense F-pure type (resp. purely F-regular type) correspond to lc (resp. plt) pairs, and Theorem \ref{F-inversion} $(2)$ suggests the following conjecture. 
\begin{conj}
Let $X$ be a non-singular variety over a field of characteristic zero and $Y=\sum_{i=1}^k t_iY_i$ a formal combination where $t_i > 0$ are real numbers and $Y_i \subsetneq X$ are closed subschemes. Let $Z \subsetneq X$ be a normal $\Q$-Gorenstein closed subvariety of codimension $r>0$ such that $Z \not\subset \cup_{i=1}^k Y_i$. 
If the pair $(X,Y+rZ)$ is plt $($resp. lc$)$ near $Z$, then the pair $(Z,Y|_Z)$ is klt $($resp. lc$)$.
\end{conj}


\begin{acknowledgement}
The author is indebted to Osamu Fujino for helpful advice about the theory of singularities of pairs, to Nobuo Hara for many discussions and some suggestions about the proof of Theorem \ref{F-inversion}, and to Kei-ichi Watanabe for useful comments. 
He is also grateful to Masayuki Kawakita and Takehiko Yasuda for valuable conversations, and to Tetsushi Ito for reading the first draft of the article.
This work was partially supported by the Japan Society for the Promotion of Science Research Fellowships for Young Scientists.
\end{acknowledgement}

\end{document}